\documentclass[12pt]{article}

\usepackage{amsmath}
\usepackage{amsfonts}
\usepackage{amsthm}

\theoremstyle{remark}

\numberwithin{equation}{section}

\begin{document}

\begin{center}
{\Large {\bf On Complete Noncompact K\"{a}hler Manifolds \ with Positive
Bisectional Curvature}} \vskip 5mm Bing-Long Chen and Xi-Ping Zhu

\smallskip 
Department of Mathematics,\\Zhongshan University,\\
Guangzhou 510275, P. R. China\\ 
\smallskip 
The Institute of Mathematical Sciences,\\The Chinese University of Hong Kong\\
Shatin, N.T., Hong Kong\\
\medskip
\end{center}

\baselineskip=20pt
\begin{abstract}
\quad We prove that a complete noncompact K\"{a}hler manifold $M^{n}$of
positive bisectional curvature satisfying suitable growth conditions is
biholomorphic to a pseudoconvex domain of {\bf C}$^{n}$ and we show that the
manifold is topologically {\bf R}$^{2n}$. In particular, when $M^{n}$ is a 
K\"{a}hler surface of positive bisectional curvature satisfying certain
natural geometric growth conditions, it is biholomorphic to {\bf C}$^{2}$.
\end{abstract}

\section*{1. Introduction}

\setcounter{section}{1} \setcounter{equation}{0} \qquad This paper is
concerned with complete noncompact K\"ahler manifolds with positive
bisectional curvature. Let us first recall some results on real Riemannian
manifolds with positive sectional curvature. Gromoll and Meyer \cite{GM} proved
that if $M$ is a complete noncompact Riemannian manifold with everywhere
positive sectional curvature, then $M$ is diffeomorphic to the Euclidean
space. Later on Greene and Wu \cite{GW} observed that this result readily follows
from the fact that a complete noncompact Riemannian manifold of positive
sectional curvature has a strictly convex exhaustion function. And the
Busemann function of Cheeger-Gromoll \cite{CG} directly gives a strictly
convex exhaustion function on a positive curved manifold. A complex manifold
is said to be Stein if it is holomorphically convex and its
global holomorphic functions separate points and give local coordinates at
every point. It is well-known that a complex Stein manifold can be
holomorphically embedded in some Euclidean space. A result of Grauert \cite{Gr}
says that a complex manifold which admits a smooth strictly plurisubharmonic
exhaustion function is Stein. Nonetheless, in the case of positive
bisectional curvature,the Busemann function does not immediately give rise
to a plurisubharmonic exhaustion function because one does not have the
geometric comparison theorem for geodesic distances as in the case of
positive sectional curvature (the Toponogov's comparison theorem). This
consideration motivated the following conjecture which was formulated by Siu
\cite{Si}.

{\bf Conjecture I:}\quad A complete noncompact K\"{a}hler manifold of positive
holomorphic bisectional curvature is a Stein manifold.

Some results concerning this conjecture were obtained. In \cite{MSY}, among other
things, Mok, Siu and Yau proved the following theorem:\vskip  3mm

{\bf \underline{Theorem ( Mok-Siu-Yau \cite{MSY} )}} \ \ \ Let $M$ be a complete
noncompact K\"ahler manifold of nonnegative bisectional curvature of complex
dimension $n\geq 2.$ If the Ricci curvatures are positive and for a fixed
base point $x_0,$there exist positive constants $C_1,C_2$ such that%
\begin{eqnarray}
(i) & Vol(B(x_0,r))\geq C_1r^{2n}\ ,& 0\leq r<+\infty\ ,\nonumber\\
(ii)& \frac{C_2^{-1}}{1+d(x_0,x)^2}\leq R(x)\leq\frac{C_2}{1+d(x_0,x)^2}\ ,
& x\in M\ ,\nonumber
\end{eqnarray}
where $B(x_0,r)$ denotes the geodesic ball of radius $r$ and centered at $%
x_0 $ , $R(x)$ denotes the scalar curvature and $d(x_0,x)$ denotes the
distance between $x_0$\ and $x$. Then $M$ is a Stein manifold. \vskip  3mm

The method used in Mok-Siu-Yau's paper is the study of the Poincar\'e-Lelong
equation on complete noncompact K\"ahler manifolds. Their result was
improved by Mok \cite{Mo2}, where the two-side-bound assumption $(ii)$\ was
replaced by the one-side-bound assumption:%
$$
(ii)^{^{\prime }}\qquad R(x)\leq \ \frac{C_2}{1+d(x_0,x)^2}\ ,\qquad x\in M\
. 
$$

In this paper we further improve their result in the following theorem:%
\vskip  3mm

{\bf \underline{Theorem 1.1}} \ \ \ Let $M$ be a complete noncompact
K\"ahler manifold of nonnegative holomorphic bisectional curvature of complex
dimension $n\geq 2$. Suppose for a fixed base point $x_0$\ there exist
positive constants $C_1,C_2$\ and $0<\varepsilon<1$\ such that%
\begin{eqnarray}
(i) & Vol(B(x_0,r))\geq C_1r^{2n}\ ,& 0\leq r<+\infty\ ,\nonumber\\
(ii)& R(x)\leq\frac{C_2}{1+d(x_0,x)^{1+\varepsilon}}\ ,& x\in M\ .\nonumber
\end{eqnarray}
Then $M$ is a Stein manifold.\vskip  3mm

Our method is the study of the following Ricci flow equation on $M$ :%
$$
\frac \partial {\partial t}g_{ij}(t)=-2R_{ij}(t)\ , 
$$
where $g_{ij}(t)$\ is a family of metrics, and $R_{ij}(t)$\ denotes the
Ricci curvature of $g_{ij}(t)$. In \cite{Sh4}, Shi established some dedicate decay
estimates for the volume element and the curvatures of the evolving metric $%
g_{ij}(t)$\ . Based on Shi's estimates we will show the injectivity radius
of the evolving metric $g_{ij}(t)$ is getting bigger and bigger and any
geodesic ball with radius less than half of the injectivity radius is almost
pseudoconvex. By a perturbation argument , we will be able to construct a
sequence of pseudoconvex domains of $M$ such that arbitrary two of these
domains form a Runge pair . Then by appealing a theorem of Markoe \cite{Ma} (see
also Siu \cite{Si}), we can deduce that $M$ is Stein.

The Gromoll and Meyer's theorem is a uniformization theorem in Riemannian
geometry category. The analogue in K\"ahler geometry is the following
well-known conjecture:

{\bf Conjecture II:}\quad A complete noncompact K\"ahler manifold of positive
holomorphic bisectional curvature of complex dimension $n$ is biholomorphic
to {\bf C}$^n$.

The first result concerning this conjecture is the following isometrically
embedding theorem of Mok, Siu and Yau \cite{MSY} and Mok \cite{Mo2}.\vskip  3mm

{\bf \underline{Theorem ( Mok-Siu-Yau \cite{MSY}, Mok \cite{Mo2} )}} \ \ \ Let $M$ be a
complete noncompact K\"ahler manifold of nonnegative holomorphic bisectional
curvature of complex dimension $n\geq 2$ . Suppose for a fixed base point $%
x_0$ :%
\begin{eqnarray}
(i) & Vol(B(x_0,r))\geq C_1r^{2n}\ ,& 0\leq r<+\infty\ ,\nonumber\\
(ii)& R(x)\leq\frac{C_2}{1+d(x_0,x)^{2+\varepsilon}}\ ,& x\in M\ ,\nonumber
\end{eqnarray}
for some $C_1,C_2>0$ and for any arbitrarily small positive constant $\varepsilon.$
Then $M$\ is isometrically biholomorphic to {\bf C}$^n$\ with the standard flat
metric.\vskip  3mm

Also in his paper \cite{Mo2}, Mok used some algebraic geometrical techniques to
control the holomorphic functions of polynomial growth on $M$ and obtained
the following holomorphic embedding theorem.\vskip  3mm

{\bf \underline{Theorem ( Mok \cite{Mo2} )}} \ \ \ Let $M$\ be a complete
noncompact K\"ahler manifold of positive holomorphic bisectional curvature
of complex dimension $n\geq 2.$\ Suppose for a fixed base point $x_0$:%
\begin{eqnarray}
(i)\ & Vol(B(x_0,r))\geq C_1r^{2n}\ ,& 0\leq r<+\infty\ ,\nonumber\\
(ii)^{\prime} & R(x)\leq\frac{C_2}{1+d(x_0,x)^2}\ ,& x\in M\ ,\nonumber
\end{eqnarray}
for some positive constants\ $C_1$ and $C_2.$Then $M$\ is biholomorphic to an affine algebraic variety. Moreover in case
of $n=2,$ if the K\"ahler manifold $M$ is actually of positive Riemannian
sectional curvature, then $M$ is biholomorphic to {\bf C}$^2.$\vskip  3mm

This result was also improved in \cite{Mo3} and \cite{To}. Instead of the bisectional
curvature, only the positivity of the Ricci curvature was assumed there.

In \cite{Sh2}, Shi announced that under the same assumptions as in the above Mok's
theorem, the K\"ahler manifold $M$ is biholomorphic to {\bf C}$^n.$ This is
actually the main theorem of Shi's Ph.D. thesis \cite{Sh3}. Here we are unwilling
to point out that there exists a gap in the proof \cite{Sh3}. \ Shi used the Ricci
flow to construct a flat K\"ahler metric on $M.$ But the flat metric may be
incomplete. Thus one can only get the biholomorphic embedding of $M$ into a
domain of {\bf C}$^n.$

In this paper ,with the help of the resolutions of generalised Poincar\'e
conjecture and a gluing argument of Shi \cite{Sh4}, we will get the following
result.\vskip  3mm

{\bf \underline{Theorem 1.2}} \ \ \ Under the same assumptions of
Theorem1.1, we have

\ \ \ $(1)$ \ \ \ \ $M$ is diffeomorphic to {\bf R}$^{2n}$ if $n>2,\ M$ is
homeomorphic to {\bf R}$^4$ if $ n=2,$ and

\ \ \ $(2)$ \ \ \ $M$ is biholomorphic to a pseudoconvex domain in {\bf C}$^n$.%
\vskip  3mm\ 

The part (2) of this result is analogous to the main theorem of Shi in
\cite{Sh4}, where the positivity of Riemannian sectional curvature was
assumed.

The combination of part (1) of Theorem 1.2 with a theorem of Ramanujam \cite{R}
immediately gives an improvement of the above Mok's
holomorphically embeddding theorem in dimension $n=2.$\ More precisely, we
have:\vskip  3mm

{\bf \underline{Corollary 1.3}} \ \ \ Let $M$ be a complete noncompact
K\"ahler surface with positive holomorphic bisectional curvature . Suppose
for a fixed base point $x_0,$there exist constants $C_1,C_2$ suth that%
\begin{eqnarray}
(i)\ & Vol(B(x_0,r))\geq C_1r^{4}\ ,& 0\leq r<+\infty\ ,\nonumber\\
(ii)^{\prime} & R(x)\leq\frac{C_2}{1+d(x_0,x)^2}\ ,& x\in M\ .\nonumber
\end{eqnarray}
Then $M$ is biholomorphic to {\bf C}$^2.$\vskip  3mm

The isometrically embedding throrem of Mok, Siu and Yau can be interpreted
as a gap phenomenon which shows that the metrics of positive holomorphic
bisectional curvature can not be too close to the flat one.Their method is
to solve the Poincar\'e-Lelong equation on complete noncompact K\"ahler
manifolds where the maximal volume growth condition (i) was assumed. Thus it
is interesting to investigate the gap phenomena on manifolds without maximal
volume growth assumption.

Recently Ni \cite{N} got some results in this direction by improving the
argument of Mok, Siu and Yau \cite{MSY}. In section 4 of this paper , we
use Shi$^{^{\prime }}$s a priori estimate and the differential Harnack
inequality of Cao \cite{Cao} for the Ricci flow to prove a general version of
Mok-Siu-Yau$^{^{\prime }}$s gap theorem.

Finally in section 5 we will give a remark on the flat metric constructed by
Shi \cite{Sh3}.

{\bf Acknowledgement}\quad We are grateful to Professor L.F. Tam for many helpful
discussions.

\section*{2. The Ricci Flow and A Priori Estimates}

\setcounter{section}{2} \setcounter{equation}{0}\qquad Let $(M,\tilde
g_{ij}) $\ be a complete noncompact \ K\"ahler manifold of complex dimension 
$n\geq 2 $ with bounded and nonnegative holomorphic bisectional curvature .
The Ricci flow is the following evolution equation for the metric 
\begin{equation}
\label{2.1}\left\{ 
\begin{array}{ll}
\bigbreak \displaystyle\frac \partial {\partial t}g_{ij}(x,t)=-2R_{ij}(x,t)\
, & on\quad M\times [0,T]\ , \\ 
g_{ij}(x,0)=\tilde g_{ij}(x)\ , & x\in M. 
\end{array}
\right. 
\end{equation}

Suppose $\{z^1,z^2,\cdots ,z^n\}$\ is the local holomorphic coordinate
system on $M$\ and $z^k=x^k+\sqrt{-1}x^{k+n},x^k,x^{k+n}\in R,\ k=1,2,...n$.
Then $\{x^1,x^2,\cdots ,x^{2n}\}$\ is the local real coordinate system on $%
M. $\ We use $i,\ j,\ k,\ l$ to denote the indices corresponding to the real
vectors and real covectors, $\alpha ,\ \beta ,\ \gamma ,\ \delta ,\cdots $
the indices corresponding to the \ holomorphic vectors and covectors. Then
the above Ricci flow equation (\ref{2.1}) can be written in holomorphic
coordinates as 
$$
\left\{ 
\begin{array}{ll}
\bigbreak \displaystyle\frac \partial {\partial t}g_{\alpha \overline{\beta }%
}(x,t)=-R_{\alpha \overline{\beta }}(x,t)\ , & on\quad M\times [0,T]\ , \\ 
g_{\alpha \bar \beta }(x,0)=\tilde g_{\alpha \bar \beta }(x)\ , & x\in M. 
\end{array}
\right. \eqno{(2.1)^\prime} 
$$

From \cite{Sh1} we know that the Ricci flow $(2.1)$ has a maximal solution $%
g_{ij}(\cdot ,t)$\ on $[0,t_{\max })$ with $t_{\max }>0$\ and the curvature
of $g_{ij}(\cdot ,t)$\ becomes unbounded as $t$ tends to $t_{\max }$ if $%
t_{\max }<+\infty \ .$ Since the maximum principle has been well understood
for noncompact manifolds, the preserving k\"ahlerity and nonnegativity of
holomorphic bisectional curvature of the Ricci flow (\ref{2.1}) were
essentially obtained by Hamilton \cite{Ha1} and Mok \cite{Mo1} (see also Shi \cite{Sh4}).

Define a function $F(x,t)$ on $M\times [0,t_{\max })$\ as follows 
\begin{equation}
\label{2.2}F(x,t)=\log \frac{\det (g_{\alpha \bar \beta }(x,t))}{\det
(g_{\alpha \bar \beta }(x,0))}\ . 
\end{equation}
\ \ \ \ \ It can be easily obtained from (2.1) that 
\begin{equation}
\label{2.3}\frac{\partial F(x,t)}{\partial t}=-R(x,t)\ . 
\end{equation}

Since the holomorphic bisectional curvature of $g_{\alpha \bar \beta }(\cdot
,t)$ is nonnegative, it follows that $F(\cdot ,t)$ is non-increasing in $t$
and $F(\cdot ,0)=0.$ And by the equation (\ref{2.1}) we know that 
\begin{equation}
\label{2.4}g_{\alpha \bar \beta }(\cdot ,t)\leq g_{\alpha \bar \beta }(\cdot
,0)\ ,\qquad on\quad M\ . 
\end{equation}
Then by definition , we have%
\begin{eqnarray}
e^{F(x,t)}R(x,t)&=& g^{\alpha \bar \beta }(x,t)R_{\alpha \bar \beta }(x,t)\cdot \frac{\det (g_{\alpha \bar \beta }(x,t))}{\det (g_{\alpha \bar \beta }(x,0))}\nonumber\\
&\leq& g^{\alpha \bar \beta }(x,0)R_{\alpha \bar \beta }(x,t)\nonumber\\
&=& g^{\alpha \bar \beta }(x,0)(R_{\alpha \bar \beta }(x,t)-R_{\alpha \bar \beta }(x,0))+R(x,0)\nonumber\\
&=& -g^{\alpha \bar \beta }(x,0)\frac{\partial ^2F(x,t)}{\partial z^\alpha \partial \bar z^\beta }+R(x,0)\nonumber\\
&=& -\Delta _0F(x,t)+R(x,0)\ ,\nonumber
\end{eqnarray}
where $\Delta _0$\ is the Laplace operator with respect to the initial
metric $g_{ij}(\cdot ,0)$ .

Combining with (\ref{2.3}) we obtain 
\begin{equation}
\label{2.5}
e^{F(x,t)}\frac{\partial F(x,t)}{\partial t}\geq \Delta
_0F(x,t)-R(x,0)\ , 
\end{equation}
and%
\begin{equation}
\label{2.6}
\Delta _0F(x,t)\leq R(x,0)\ ,\qquad on\quad M\times [0,t_{\max })\ . 
\end{equation}

In the PDE jargon, if the scalar curvature $R(x,0)$ of the initial
metric satisfies suitable growth conditions, the differential inequalities
(2.5) and (2.6) will give two opposite estimates of $F$ by its average. Shi
\cite{Sh4} observed that the combination of these two opposite estimates will give
the following a priori estimate for the function $F.$\vskip  3mm

{\bf \underline{Lemma 2.1 (Shi \cite{Sh4})}} \ \ \ Suppose $(M,\tilde g_{\alpha
\bar \beta })$ is a complete noncompact \ K\"ahler manifold of complex
dimension $n\geq 2$ with bounded and nonnegative holomorphic bisectional
curvature. Suppose there exist positive constants $C_1,C_2$ and $0<\varepsilon<1$ such
that for any $x_0,$%
\begin{eqnarray}
(i) & R(x,0)\leq C_1\ ,& x\in M\ ,\nonumber\\
(ii)& \frac 1{Vol(B_0(x_0,r))}\int\nolimits_{B_0(x_0,r)}R(x,0)dV_0\leq \frac{C_2}{1+r^{1+\varepsilon}}\ ,
& 0\leq r<+\infty\ ,\nonumber
\end{eqnarray}
where $B_0(x_0,r)$ is the geodesic ball of radius $r$ and centered at $x_0$
with respect to the metric $\tilde g_{\alpha \bar \beta }(x).$ Then the
solution of (2.1)$^{\prime }$ satisfies the estimate:%
$$
F(x,t)\geq -C(t+1)^{\frac{1-\varepsilon}{1+\varepsilon}}\ ,\qquad on~M\times [0,t_{\max })\ , 
$$
where $0<C<+\infty $ is a constant depending only on $n,$ $\varepsilon,\ C_1$and \ $%
C_2.$\vskip  3mm

{\bf \underline{Proof.}} \ \ \ Since we will use the two opposite estimates
later in this paper, \ for convenience, we sketch the proof here.

Without loss of generality , by replacing $M$ by $M\times${\bf C}$^2$ if necessary,
we may assume that the dimension of $M$\ is $\geq 4$ and 
\begin{equation}
\label{2.7}C_3\left( \frac{r_2}{r_1}\right) ^4\leq \frac{Vol(B_0(x,r_2))}{%
Vol(B_0(x,r_1))}\leq \left( \frac{r_2}{r_1}\right) ^{2n}\ ,\qquad 0\leq
r_1\leq r_2<+\infty \ , 
\end{equation}
by the standard volume comparison, where $C_3$\ is a positive constant
depending only on $n.$

Denote $\tilde \nabla $\ to be the covariant derivatives with respect to the
initial metric $\tilde g_{\alpha \bar \beta }(x)$\ . Since the holomorphic
bisectional curvature of$\ \tilde g_{\alpha \bar \beta }(x)$\ is
nonnegative, we know that the Ricci curvature is also nonnegative. Then by
Theorem 1.4.2 of Schoen and Yau \cite{SY} and a simple scaling argument, there
exists a constant $C(n)>0$\ depending only on $n$\ such that for any fixed
point $x_0\in M$ and any number $0<a<+\infty ,$\ there exists a smooth
function $\varphi (x)\in C^\infty (M)$\ satisfying 
\begin{equation}
\label{2.8}\left\{ 
\begin{array}{l}
\bigbreak \displaystyle e^{-C(n)\left( 1+\frac{d_0(x,x_0)}a\right) }\leq
\varphi (x)\leq \ \ e^{-\left( 1+\frac{d_0(x,x_0)}a\right) }\ , \\ \bigbreak %
\left| \tilde \nabla \varphi (x)\right| _0\leq \frac{C(n)}a\varphi (x)\ , \\ 
\left| \Delta _0\varphi (x)\right| _0\leq \frac{C(n)}{a^2}\varphi (x)\ , 
\end{array}
\right. 
\end{equation}
for $x\in M,$\ where $d_0(x,x_0)$\ is the distance between $x$\ and $x_0$\
with respect to the metric $\tilde g_{\alpha \bar \beta }(x)$\ and $\left|
\cdot \right| _0$ is the corresponding norm.

By the volume growth (\ref{2.7}), the manifold $(M,\tilde g_{\alpha \bar
\beta })$ is parabolic. Let $G_0(x,y)$ be the positive Green's
function on $(M,$ $\tilde g_{\alpha \bar \beta })$\ . From Li-Yau's estimate
for Green's function and (\ref{2.7})\ we know that for $x,y\in M,$%
\begin{equation}
\label{2.9}\left\{ 
\begin{array}{l}
\bigbreak \displaystyle\frac{C_4^{-1}d_0(x,y)^2}{Vol(B_0(x,d_0(x,y)))}\leq
G_0(x,y)\leq\frac{C_4d_0(x,y)^2}{Vol(B_0(x,d_0(x,y)))}\ , \\ \left|
\tilde \nabla G_0(x,y)\right| _0\leq\frac{C_5d_0(x,y)}{%
Vol(B_0(x,d_0(x,y)))}\ , 
\end{array}
\right. 
\end{equation}
for some positive constants $C_4,C_5$ depending only on $n.$

Denote%
$$
\Omega _\alpha =\{y\in M|G_0(x,y)>\alpha \}\qquad for\ \ any\quad \alpha
>0\ . 
$$

By the differential inequality (\ref{2.6}), we have for any $%
x_0\in M$ and any $t\in [0,t_{\max }),$%
$$
F(x_0,t)=\int\nolimits_{\Omega _\alpha }\left( \alpha -G_0(x_0,y)\right) \Delta _0F(y,t)dV_0(y)-\int\nolimits_{\partial \Omega _\alpha }F(y,t)\frac{\partial G_0(x_0,y)}{\partial\nu}d\sigma (y)
$$
\begin{equation}
\label{2.10}
\qquad\quad\geq -\int\nolimits_{\Omega _\alpha }G_0(x_0,y)R(y,0)dV_0(y)-\int\nolimits_{\partial \Omega _\alpha }F(y,t)\frac{\partial G_0(x_0,y)}{\partial\nu}d\sigma (y)\ .
\end{equation}
Let $r(\alpha )\geq 1$ be a number such that 
\begin{equation}
\label{2.11}\frac{r(\alpha )^2}{Vol(B_0(x_0,r(\alpha )))}=\alpha \ . 
\end{equation}
Then by (\ref{2.7}), (\ref{2.9}) and the assumption (ii) , we
have 
\begin{equation}
\label{2.12}C_6^{-1}r(\alpha )\leq d_0(x_0,y)\leq C_6r(\alpha )\ ,\qquad
for\quad any\quad y\in \partial \Omega _\alpha \ , 
\end{equation}
\begin{equation}
\label{2.13}-\int\nolimits_{\partial \Omega _\alpha }F(y,t)\frac{\partial
G_0(x_0,y)}{\partial\nu}d\sigma (y)\geq \frac{C_5r(\alpha )}{%
Vol(B_0(x_0,r(\alpha )))}\ \int\nolimits_{\partial \Omega _\alpha
}F(y,t)d\sigma (y)\ , 
\end{equation}
and%
\begin{eqnarray}
\label{2.14}
-\int\nolimits_{\Omega _\alpha }G_0(x_0,y)R(y,0)dV_0(y) &\geq& -\int\nolimits_{B(x_0,C_6r(\alpha ))}G_0(x_0,y)R(y,0)dV_0(y)\nonumber\\
&\geq& -C_7(r(\alpha ))^{1-\varepsilon}\ ,
\end{eqnarray}
where $C_6,C_7$ are positive constants depending only on $C_1,C_2$ and $n$ .

Substituting (\ref{2.13}), (\ref{2.14}) into (\ref{2.10}) and integrating
from $\frac \alpha 2$ to $\alpha $ , one readily get 
$$
F(x_0,t)\geq -C_8r(\alpha )^{1-\varepsilon}+\frac{C_8}{Vol(B_0(x_0,r(%
\alpha )))}\ \int\nolimits_{B_0(x_0,C_6r(\alpha ))}F(y,t)dV_0(y)\ , 
$$
for some positive constant depending only on $\varepsilon,C_1,C_2$\ and $n$ .

Then by (2.12) and (2.7) , we obtain that for any $a>0,$%
\begin{equation}
\label{2.15}
F(x_0,t)\geq -C_9a^{1-\varepsilon}+\frac{C_9}{Vol(B_0(x_0,a))}\int\nolimits_{B_0(x_0,a)}F(y,t)dV_0(y)\ , 
\end{equation}
where $C_9$ is a positive constant depending only on $\varepsilon,C_1,C_2$ and $n$.

On the other hand , by multiplying the differential inequality (\ref{2.5})
by the function $\varphi (x)$ in (\ref{2.8}) and integrating by parts, we get%
\begin{eqnarray}
\frac \partial {\partial t}\int\nolimits_M\varphi (x)e^{F(x,t)}dV_0 &\geq&\int\nolimits_M(\Delta _0F(x,t)-R(x,0))\varphi (x)dV_0\nonumber\\
&\geq&\frac{C(n)}{a^2}\int\nolimits_MF(x,t)\varphi (x)dV_0-\int\nolimits_MR(x,0)\varphi (x)dV_0\ .\nonumber
\end{eqnarray}
By (2.8) and the assumption (i) (ii), it is easy to see that%
$$
\int\nolimits_MR(x,0)\varphi (x)dV_0\leq \frac{C_{11}}{a^{1+\varepsilon}}%
Vol(B_0(x_0,a)) 
$$
for some positive constant depending only on $n$ . Then we have 
\begin{equation}
\label{2.16}\int\nolimits_M\varphi (x)(1-e^{F(x,t)})dV_0\leq \frac{C_{11}t}{%
a^{1+\varepsilon}}Vol(B_0(x_0,a))+\frac{C(n)t}{a^2}\int\nolimits_M(-F(x,t))\varphi
(x)dV_0. 
\end{equation}

Denote%
$$
F_{\min }(t)=\inf \{F(x,t)|x\in M\}\ . 
$$
It is clear that 
\begin{equation}
\label{2.17}\int\nolimits_M\varphi (x)(1-e^{F(x,t)})dV_0\geq \frac
1{2(1-F_{\min }(t))}\int\nolimits_M(-F(x,t))\varphi (x)dV_0\ . 
\end{equation}
Therefore we get from (\ref{2.16}) and (\ref{2.17}), 
\begin{equation}
\label{2.18}\frac 1{Vol(B_0(x_0,a))}\int\nolimits_{B_0(x_0,a)}F(x,t)dV_0\geq
-C_{12}t(1-F_{\min }(t))\left( \frac 1{a^{1+\varepsilon}}-\frac{F_{\min }(t)}{a^2}%
\right) \ , 
\end{equation}
for any $a>0$ and $t\in [0,t_{\max }),$ where $C_{12}$ is a positive
constant depending only on $\varepsilon,C_1,C_2$ and $n$.

Combining (\ref{2.15}) and (\ref{2.18}) , we deduce that 
\begin{equation}
\label{2.19}F_{\min }(t)\geq -C_9a^{1-\varepsilon}-C_9C_{12}(1-F_{\min }(t))\left(
\frac t{a^{1+\varepsilon}}-\frac t{a^2}F_{\min }(t)\right) \ , 
\end{equation}
for any $a>0$ and $t\in [0,t_{\max })$ . By taking $a=C_{13}(t+1)^{\frac
12}(-F_{\min }(t))^{\frac 12}$ with $C_{13}$ large enough , we then get the
desired estimate%
$$
F_{\min }(t)\geq -C(t+1)^{\frac{1-\varepsilon}{1+\varepsilon}}\ ,\qquad for~t\in [0,t_{\max })\
. 
$$
\hfill
$Q. E. C. $

Next we are going to prove that the maximal volume growth condition is
preserved under the Ricci flow (\ref{2.1}) (or (2.1)$^{^{\prime }}$) . More
precisely, we have\vskip  3mm

{\bf \underline{Lemma 2.2}} \ \ \ Suppose $(M,\tilde g_{\alpha \bar \beta })$
is a complete noncompact K\"ahler manifold of complex dimension $n\geq 2$
with bounded and nonnegative bisectional curvature. Suppose for a fixed base
point $x_0$ there exist positive \ constant $C_1,C_2$ and $0<\varepsilon<1$ such that%
\begin{eqnarray}
(i) & Vol(B(x_0,r))\geq C_1r^{2n}\ ,& 0\leq r<+\infty\ ,\nonumber\\
(ii)& R(x)\leq\frac{C_2}{1+d(x_0,x)^{1+\varepsilon}}\ ,& x\in M\ .\nonumber
\end{eqnarray}
Let $g_{\alpha \bar \beta }(\cdot ,t)$ be the solution of the Ricci flow
(2.1)$^{^{\prime }}$ with $\tilde g_{\alpha \bar \beta }$ as initial metric,
and let $Vol_t(B_t(x_0,r))$ be the volume of the geodesic ball of radius $r$
and centered at $x_0$ with respect to the metric $g_{\alpha \bar \beta
}(\cdot ,t)$. Then 
\begin{equation}
\label{2.20}Vol_t(B_t(x_0,r))\geq C_1r^{2n} 
\end{equation}
for all $t\in [0,t_{\max })$ and $0\leq r<+\infty .$\vskip  3mm

{\bf \underline{Proof.}} \ \ \ We have already noticed that the Ricci
curvature of the solution $g_{\alpha \bar \beta }(\cdot ,t)$ is nonnegative
. That says the metric is shrinking under the flow (2.1)$^{^{\prime }}.$
Thus we have%
\begin{eqnarray}
\label{2.21}
Vol_t(B_t(x_0,r)) &\geq&Vol_t(B_0(x_0,r))\nonumber\\
&=&\int\nolimits_{B_0(x_0,r)}e^{F(x,t)}dV_0\nonumber\\
&=&Vol(B_0(x_0,r))+\int\nolimits_{B_0(x_0,r)}(e^{F(x,t)}-1)dV_0
\end{eqnarray}
Clearly the assumptions of Lemma 2.2 are stronger then those of Lemma 2.1.
Then we can use (2.8) and (2.16) to deduce%
\begin{eqnarray}
\int\nolimits_{B_0(x_0,r)}(e^{F(x,t)}-1)dV_0&\geq&C_3\int\nolimits_M(e^{F(x,t)}-1)\varphi (x)dV_0\nonumber\\
&\geq&C_4t\left[\frac{F_{\min }(t)}{r^2}-\frac 1{r^{1+\varepsilon}}\right]\cdot Vol(B_0(x_0,r))\nonumber
\end{eqnarray}for some positive constants $C_3,C_4$\ depending on $C_1,C_2$\
and $n$.

Substituting the above inequality into the right hand side of (\ref{2.21})
and dividing by $r^{2n}$, we have%
$$
\lim \limits_{r\rightarrow +\infty }\frac{Vol_t(B_t(x_0,r))}{r^{2n}}\geq \ \
\lim \limits_{r\rightarrow +\infty }\frac{Vol(B_0(x_0,r))}{r^{2n}}\geq C_1\
. 
$$
Hence the standard volume comparison theorem implies (\ref{2.20}).\hfill
$Q. E. C. $

\section*{3. The Proof of Theorem 1.1 and 1.2}

\setcounter{section}{3} \setcounter{equation}{0}\ ~\quad In this section we
will prove Theorem 1.1 and 1.2 simutaneously.

First, let us recall the local injectivity radius estimate of Cheeger,
Gromov and Taylor \cite{CGT}, which says that for any complete Riemannian manifold 
$N$ of dimension $m$ with $\lambda \leq $ sectional curvatures of $N\leq\Lambda,$
let $r$ be a positive constant and $r<\frac \pi {4\sqrt{\Lambda }}$
if $\Lambda >0,$ the injetivity radius of $N$ at a point $P$ can be bounded
from below as follows 
\begin{equation}
\label{3.1}\ inj_N(P)\geq r\frac{Vol(B(P,r))}{Vol(B(P,r))+V^{2m}(2r)}\ , 
\end{equation}
where $V^m(2r)$ denotes the volume of a ball of radius $2r$ in the $m-$%
dimensional model space $V^m$ with constant sectional curvature $\lambda .$

In particular , it implies that for a complete Riemannian manifold $N$ of
dimension $m$ with the sectional curvature bounded between $-1$ and $1$ ,
the injectivity radius at a point $P$ can be estimated as follows 
\begin{equation}
\label{3.2}\ inj_N(P)\geq \frac 12\frac{Vol(B(P,\frac 12))}{Vol(B(P,\frac
12))+V} 
\end{equation}
for some positive constant $V$ depending only on $m.$ Further, if in
addition $N$ satisfies the maximal volume growth condition%
$$
Vol(B(x_0,r))\geq C_1r^m\ ,\qquad 0\leq r\leq +\infty \ , 
$$
then (\ref{3.2}) gives 
\begin{equation}
\label{3.3}inj_N(P)\geq C>0 
\end{equation}
for some positive constant $C$ depending only on $C_1$ and $m$.

Consider $(M,\tilde g_{\alpha \bar \beta })$ to be a complete noncompact
K\"ahler manifold with complex dimension $n\geq 2$ and satisfying the
assumptions of Theorem 1.1 (or Theorem 1.2). \ Let $g_{\alpha \bar \beta
}(\cdot ,t)$ be the maximal solution of the Ricci flow (2.1)$^{^{\prime }}$.
Lemma 2.1 tells us that 
$$
\frac{{\rm det}(g_{\alpha\bar{\beta}}(x, t))}{{\rm det}(g_{\alpha\bar{\beta}}(x, 0))}
\geq e^{-C(t+1)^{\frac{1-\varepsilon}{1+\varepsilon}}},
\ \ \ \mbox{\rm on }
M\times[0, t_{\rm max}),
$$
which together with (2.4) implies
$$
g_{\alpha\bar{\beta}}(x, 0)\geq
g_{\alpha\bar{\beta}}(x, t)\geq
e^{-C(t+1)^{\frac{1-\varepsilon}{1+\varepsilon}}}g_{\alpha\bar{\beta}}(x, 0),
\ \ \ \mbox{\rm on }
M\times[0, t_{\rm max}).
$$
Since the Ricci flow equation (2.1)$^{^{\prime }}$ is the parabolic version of the complex Monge-Amp\`{e}re eqution on the K\"{a}hler manifold, the above inequality is corresponding to the second order estimate for the Monge-Amp\`{e}re eqution. It is well known that the third 
order and higher order estimates for Monge-Amp\`{e}re were developed by Calabi and Yau.
Similarly, by adapting the Calabi and Yau's arguments, Shi proved in
\cite{Sh4} that the derivative and higher order estimates for $
g_{\alpha\bar{\beta}}(x, t)$ are uniformly bounded on any finite time interval. 
In particular, this implies that the solution  $
g_{\alpha\bar{\beta}}(\cdot, t)$ exists for all $t\in [0, +\infty)$. On the other hand,
by applying the differential Harnack inequality of Cao \cite{Cao}, we know that $%
tR(\cdot,t)$ is nondecreasing in time. It then follows from (\ref{2.3}) that
for $x\in M$ and $t\in [0,+\infty)$, %
\begin{eqnarray}
\label{3.4}
-F(x,2t)&=&\int\nolimits_0^{2t}R(x,s)ds\nonumber\\
&\geq&\int\nolimits_t^{2t}R(x,s)ds\nonumber\\
&\geq&tR(x,t)\int\nolimits_t^{2t}\frac 1sds\nonumber\\
&=&(\log 2)\cdot tR(x,t)\ .
\end{eqnarray}Combining (\ref{3.4}) with Lemma 2.1, we get 
\begin{equation}
\label{3.5}R(x,t)\leq C_3(t+1)^{\frac{-2\varepsilon}{1+\varepsilon}}\ , 
\end{equation}
for some positive constant $C_3$ depending only on $\varepsilon,C_1,C_2$ and $n.$

Moreover by using the derivative estimates of Shi (see \cite{Sh1} or Theorem 7.1
of Hamilton \cite{Ha1}), we have 
\begin{equation}
\label{3.6}\left| \nabla ^pR_{ijkl}(x,t)\right| \leq
C(\varepsilon,C_1,C_2,n,p)(t+1)^{- \frac{\varepsilon}{1+\varepsilon}(p+2)}\ , 
\end{equation}
for $x\in M,t\geq 1$ and any\ integer $p\geq 0.$

We have shown in Lemma 2.2 that the maximal volume growth condition is
preserved under the flow (2.1)$^{^{\prime }}.$ By a standard scaling
argument , it follows from (\ref{3.3}) and (\ref{3.5}) that the injectivity
radius of $M$ with respect to the metric $g_{\alpha \bar \beta }(\cdot ,t)$
has the following estimate 
\begin{equation}
\label{3.7}inj(M,g_{\alpha \bar \beta }(\cdot,t))\geq C_4(t+1)^{\frac{\varepsilon}{1+\varepsilon}}\
, 
\end{equation}
where $C_4$ is a positive constant depending only on $\varepsilon,\ C_1,\ C_2$ and $n$.

Since the Ricci curvature of $g_{\alpha \bar \beta }(\cdot ,t)$ is
nonnegative for all $x\in M$ and $t\geq 0,$ we know from the equation (2.1)$%
^{^{\prime }}$ that the ball $B_t\left(x_0,\frac{C_4}2(t+1)^{\frac{\varepsilon}{1+\varepsilon}}\right),$ of
radius $\frac{C_4}2(t+1)^{\frac{\varepsilon}{1+\varepsilon}}$ with respect to the metric $%
g_{\alpha \bar \beta }(\cdot,t),$ contains the ball $B_0\left(x_0,\frac{C_4}2(t+1)^{\frac{\varepsilon}{1+\varepsilon}}\right),$ of the same radius with
respect to the metric\ $\tilde g_{\alpha \bar \beta }=g_{\alpha \bar \beta
}(\cdot ,0)$.Thus we deduce from (\ref{3.7}) that 
\begin{equation}
\label{3.8}\pi _p(M)=0\qquad and\qquad \pi _q(M,\infty )=0 
\end{equation}
for any $p\geq 1,1\leq q\leq 2n-2,$ where $\pi _q(M,\infty )$ is the $q-$th
homotopy group of $M$ at infinity.

Then by the resolutions of generalised Poincar\'{e} conjecture ( see \cite{F}, 
\cite{Sm} ), we know that $M$ is homeomorphic to {\bf R}$^{2n}.$
Notice Gompf$^{^{\prime }}$s result says that among the Euclidean spaces
only {\bf R}$^{4}$ has exotic differential structures. So for $n>2$ the
homeomorphisms can be maken to be diffeomorphisms. This gives the proof of
the part (1) of Theorem 1.2 .

Also the injectivity radius estimate (\ref{3.7}) tells us the exponential
map provides a diffeomorphism between the balls of $M$ and the Euclidean
space. In the following we want to modify the exponential maps to become
biholomorphims.

Let $T_{x_0}M$ denote the real tangent space of $M$ at $x_0$, $J_M$ denote
the complex structure of $T_{x_0}M$. For fixed $t,$ choose a standard
orthonormal basis $\{e_1,J_Me_1,\cdot \cdot \cdot ,e_n,J_Me_n\}$ of $%
T_{x_0}M $ . Since the metric $g_{\alpha \bar \beta }(\cdot ,t)$ is a family
of smooth K\"ahler metrics, we may assume the basis is smooth in time.

For any $v\in T_{x_0}M$ ,one can write%
$$
v=x_1e_1+y_1J_Me_1+...+x_ne_n+y_nJ_Me_n\ . 
$$
We now construct a real linear isomorphism $L:T_{x_0}M\rightarrow${\bf C}$^n$
defined by%
$$
L(v)=(z_1,z_2,...z_n)\in {\bf C}^n\ , 
$$
where $z_i=x_i+\sqrt{-1}y_i$ , $i=1,2,\cdots ,n.$\ It is also clear that $L$
varies smoothly in $t$. Equip {\bf C}$^n$ with the standard flat K\"ahler metric.
Let $\hat \nabla $ denote the covariant derivative and  $\hat B(0,r)$ denote
the ball of radius $r$ with respect to this standard metric.

Let us use $\exp {}_{x_0}^t$ to denote the exponential map with respect to
the metric $g_{\alpha \bar \beta }(\cdot ,t)$. By (\ref{3.7}), the map $%
\varphi _t=\exp {}_{x_0}^t\circ L^{-1}$ is a diffeomorphism from $\hat B\left(0, 
\frac{C_4}2(t+1)^{\frac{\varepsilon}{1+\varepsilon}}\right)$ to the geodesic ball $B_t\left(x_0,\frac{C_4}%
2(t+1)^{\frac{\varepsilon}{1+\varepsilon}}\right)$ of $M$ with repect to the metric $g_{\alpha \bar
\beta }(\cdot ,t)$, and the map is nonsingular on $\hat B\left(0,\frac{C_4}%
2(t+1)^{\frac{\varepsilon}{1+\varepsilon}}\right)$ .Let us write the solution $g_{\alpha \bar \beta
}(\cdot ,t)$ in real coordinates as $g_{ij}(\cdot ,t)$ . We consider the
pull back metric 
\begin{equation}
\label{3.9}\varphi _t^{*}(g_{ij}(\cdot ,t))=g_{ij}^{*}(\cdot ,t)dx^idx^j\
,\qquad on\quad \hat B\left(0,\frac{C_4}2(t+1)^{\frac{\varepsilon}{1+\varepsilon}}\right)\ . 
\end{equation}
Clearly we can also write the pull back metric in the complex coordinates of 
{\bf C}$^n$ as 
\begin{equation}
\label{3.10}\varphi _t^{*}(g_{ij}(\cdot ,t))=g_{AB}^{*}(\cdot ,t)dz^Adz^B\
,\qquad on\quad \hat B(0,\frac{C_4}2(t+1)^{\frac{\varepsilon}{1+\varepsilon}})\ , 
\end{equation}
where $A,B=\alpha $ or $\bar \alpha $($\alpha =1,2,...n).$ Since $\varphi _t$
is not holomorphic in general , the metric $g_{AB}^{*}(\cdot ,t)$ is not
Hermitian with respect to the standard complex structure of {\bf C}$^n$

Notice that $g_{AB}^{*}(\cdot ,t)$ is just the representation of the metric $%
g_{\alpha \bar \beta }(\cdot ,t)$ in geodesic coordinates .The following
lemma due to Hamilton ( Theorem 4.10 in \cite{Ha2} ) is useful to estimate $%
g_{AB}^{*}(\cdot ,t)$.\vskip  3mm

{\bf \underline{Lemma 3.1}} \ \ \ Suppose the metric $g_{ij}dx^idx^j$ is in
geodesic coordinates. Suppose the Riemannian curvature $R_m$ is bounded
between $-B_0$ and $B_0$. Then there exist positive constants $c,C_0$
depending only on the dimension such that for any $\left| x\right| \leq
\frac c{\sqrt{B_0}}$, the following holds%
$$
\left| g_{ij}-\delta _{ij}\right| \leq C_0B_0\left| x\right| ^2\ . 
$$
Furthermore, if in addition $\left| \nabla R_m\right| \leq B_0$ and $\left|
\nabla ^2R_m\right| \leq B_0$ , then%
$$
\left| \frac \partial {\partial x^j}g_{kl}\right| \leq C_0B_0\left| x\right|
\qquad and\qquad \left| \frac{\partial ^2}{\partial x^i\partial x^j}%
g_{kl}\right| \leq C_0B_0 
$$
for any $\left| x\right| \leq \frac c{\sqrt{B_0}}.$

Applying this lemma to the evolving metric $g_{\alpha \bar \beta }(\cdot ,t)$
, we can find positive constants $\bar c,C_5$ depending only on $\varepsilon,C_1,C_2$
and $n$ such that for any $z\in \hat B\left(0,\bar c(t+1)^{\frac{\varepsilon}{1+\varepsilon}}\right)$ and $%
t\in [1,+\infty )$,%
\begin{eqnarray}
\label{3.11}\left| g_{\alpha \bar \beta }^{*}(z,t)-\delta _{\alpha \beta }\right| _t&\leq&C_5\left| z\right| ^2(t+1)^{-\frac{2\varepsilon}{1+\varepsilon}}\ ,\\
\label{3.12}\left| g_{\bar \alpha \beta }^{*}(z,t)-\delta _{\alpha \beta }\right| _t&\leq&C_5\left| z\right| ^2(t+1)^{-\frac{2\varepsilon}{1+\varepsilon}}\ ,\\
\label{3.13}\left| g_{\alpha \beta }^{*}(z,t)\right| _t\quad&\leq&C_5\left| z\right| ^2(t+1)^{-\frac{2\varepsilon}{1+\varepsilon}}\ ,\\
\label{3.14}\left| g_{\bar \alpha \bar \beta}^{*}(z,t)\right| _t\quad&\leq&C_5\left| z\right| ^2(t+1)^{-\frac{2\varepsilon}{1+\varepsilon}}\ ,\\
\label{3.15}\left| \hat \nabla g_{AB}^{*}(z,t)\right| _t\ \ &\leq&C_5\left| z\right|(t+1)^{-\frac{2\varepsilon}{1+\varepsilon}}\ ,\\
\label{3.16}\left| \hat \nabla\hat \nabla g_{AB}^{*}(z,t)\right| _t&\leq&C_5(t+1)^{-\frac{2\varepsilon}{1+\varepsilon}}\ ,
\end{eqnarray}where $\left|\cdot\right| _t$ is the norm with respect to the
metric $g_{ij}^{*}(z,t).$

By a further restriction on the small constant $\bar c,$ we may assume that
in the real coordinates 
\begin{equation}
\label{3.17}\frac 12\ g_{ij}^{*}(z,t)\leq \delta _{ij}\leq 2g_{ij}^{*}(z,t)\
. 
\end{equation}

Let $\varphi _{t}^{\ast }J_{M}$\ and $\bar{\partial}^{t}=\varphi _{t}^{\ast
}(\bar{\partial})$ be the pull back complex structure and $\bar{\partial}$%
-operator of $M$ on $\hat{B}\left(0,\bar{c}(t+1)^{\frac{\varepsilon }{%
1+\varepsilon }}\right)$. It is obvious that $\varphi _{t}$ is holomorphic with
respect to the pull back complex structure $\varphi _{t}^{\ast }J_{M}$\ .But
the functions $z^{\alpha }$ $(\alpha =1,2,...n)$ are not holomorphic with
respect to $\varphi _{t}^{\ast }J_{M}$\ in general. This just indicates the
difference of the pull back complex structure $\varphi _{t}^{\ast }J_{M}$\ \
with the standard complex srtucture $J_{{\bf C}^{n}}$ on $\hat{B}\left(0,\bar{c}(t+1)^{
\frac{\varepsilon }{1+\varepsilon }}\right)$. However we can estimate the
difference as follows.

Let us denote $\{x^1,...,x^n,x^{n+1},...x^{2n}\}=\{x^1,...,x^n,y^1,...y^n\}$
as real coordinates for $\hat B(0,\bar c(t+1)^{\frac \varepsilon{1+\varepsilon}})$ . We then
write%
$$
\varphi _t^{*}J_M=J_j^i(x,t)\frac \partial {\partial x^i}\otimes dx^j\qquad
and\qquad J_{{\bf C}^n}=\hat J_j^i\frac \partial {\partial x^i}\otimes dx^j\ . 
$$
By definition,$J_j^i(x,t)$ is just the representation of the complex
structure $J_M$ in the normal coordinate at $x_0,$and $J_j^i(0,t)=\hat
J_j^i(0).$ Denote $\nabla ^t$ and $\Gamma _{ij}^{t_k}$ by the covariant
derivative and Christoffel symboles with respect to the pull back metric $%
g_{ij}^{*}(\cdot ,t)$of \ $g_{\alpha \bar \beta }(\cdot ,t)$ on $\hat
B\left(0,\bar c(t+1)^{\frac \varepsilon{1+\varepsilon}}\right)$.

Set%
$$
H_k^j=\hat J_k^j-J_k^j\ . 
$$
By the kahlerity of $J_M,$ we have 
\begin{equation}
\label{3.18}x^i\nabla _i^tH_k^j=x^i\Gamma _{ip}^{t_j}\hat J_k^p-x^i\Gamma
_{ik}^{t_p}\hat J_p^j\ . 
\end{equation}

Since the metric $g_{ij}^{*}(\cdot ,t)$ is actually the representation of $%
g_{\alpha \bar \beta }(\cdot ,t)$ in geodesic coordinates, it follows from
the Gauss Lemma (see also Lemma 4.1 of \cite{Ha2}) that 
\begin{equation}
\label{3.19}g_{ij}^{*}x^i=\delta _{ij}x^i\ ,\qquad on\quad \hat B(0,\bar
c(t+1)^{\frac{\varepsilon}{1+\varepsilon}})\ . 
\end{equation}
As in \cite{Ha2}, we introduce the symmetric tensor%
$$
A_{ij}=\frac 12x^k\frac \partial {\partial x^k}g_{ij}^{*}\ . 
$$

By (\ref{3.19}), we get 
$$
x^k\frac \partial {\partial x^i}g_{jk}^{*}=\delta
_{ij}-g_{ij}^{*}=x^k\frac \partial {\partial x^j}g_{ik}^{*}\ , 
$$
and hence from the formula for $\Gamma _{ij}^{t_k}$, it follows%
\begin{equation}
\label{3.20}x^j\Gamma _{jk}^{t_i}=g^{*il}A_{kl}\ . 
\end{equation}

Combining (\ref{3.17}), (\ref{3.18}) and (\ref{3.20}), we get 
\begin{equation}
\label{3.21}\left| x^i\nabla _i^tH_k^j\right| _t\leq C_6\left| A_{pq}\right|
_t\ ,\qquad on\quad \hat B\left(0,\bar c(t+1)^{\frac{\varepsilon}{1+\varepsilon}}\right)\ , 
\end{equation}
where $C_6$ is a positive constant depending only on $n.$

The following lemma due to Hamilton \cite{Ha2} gives an estimate for $\left|
A_{pq}\right| $.\vskip  3mm

{\bf \underline{Lemma 3.2}} \ \ \ There exist constants $c>0$ and $%
C_0<\infty $ such that if the metric $g_{ij}$ is in geodesic coordinates with 
$\left| R_m\right| \leq B_0$ in the ball of radius $r\leq \frac c{\sqrt{B_0}%
}$, then%
$$
\left| A_{ij}\right| \leq C_0B_0r^2\ . 
$$
\vskip  3mm

{\bf \underline{Proof.}} \ \ \ This is Theorem 4.5 of \cite{Ha2}.\hfill
$Q. E. C. $

In our context, this lemma implies 
\begin{equation}
\label{3.22}\left| x^i\nabla _i^tH_k^j\right| _t\leq C_7\left| x\right|
^2(t+1)^{-\frac{2\varepsilon}{1+\varepsilon}}\ ,\qquad on\quad \hat B(0,\bar c(t+1)^{\frac{\varepsilon}{1+\varepsilon%
}})\ , 
\end{equation}
for some positive contant $C_7$ depending only on $\varepsilon,C_1,C_2$ and $n$.

Set%
$$
M(r)=\sup \limits_{\{\left| x\right| \leq r\}}\left| H_k^j\right| _t\ . 
$$

Suppose in the ball $\left| x\right| \leq r,$ $\left| H_k^j\right| _t$
achieves its maximum at $x_0.$ Then by (3.22) and the fact that $H_k^j(0)=0,$
we have%
\begin{eqnarray}
M(r)^2&=&\int\nolimits_0^1\frac \partial {\partial s}\left| H_k^j(sx_0^1,\cdot \cdot \cdot ,sx_0^{2n})\right| _t^2ds\nonumber\\
&=&2\int\nolimits_0^1\langle H_k^j(sx_0^1,\cdot \cdot \cdot ,sx_0^{2n}),x_0^i\nabla _i^tH_k^j(sx_0^1,\cdot \cdot \cdot ,sx_0^{2n})\rangle _tds\nonumber\\
&\leq&M(r)\cdot C_7r^2(t+1)^{-\frac{2\varepsilon}{1+\varepsilon}}\ .\nonumber
\end{eqnarray}Therefore we get the following estimate for the difference of
the two complex structures $\varphi _t^{*}J_M$ and $J_{{\bf C}^n},$%
\begin{equation}
\label{3.23}M(r)\leq C_7r^2(t+1)^{-\frac{2\varepsilon}{1+\varepsilon}}\qquad on\quad \hat
B(0,\bar c(t+1)^{\frac{\varepsilon}{1+\varepsilon}})\ . 
\end{equation}
Since $z^\alpha (\alpha =1,2,...n)$ are holomorphic with respect to $%
J_{{\bf C}^n}, $ we get from (\ref{3.17}) and (\ref{3.23}) that%
\begin{eqnarray}
\label{3.24}
\left| \bar \partial ^tz^\alpha \right| _t&\leq&\sum\limits_i\left| \frac \partial {\partial x^i}z^\alpha +\sqrt{-1}\varphi _t^{*}J_M(\frac \partial {\partial x^i})z^\alpha \right| _t\nonumber\\
&\leq&\sum\limits_i\left| \frac \partial {\partial x^i}z^\alpha +\sqrt{-1}J_{{\bf C}^n}(\frac \partial {\partial x^i})z^\alpha \right| _t+\sum\limits_i\left| (\varphi _t^{*}J_M-J_{{\bf C}^n})(\frac \partial {\partial x^i})z^\alpha \right|_t \nonumber\\
&=&\sum\limits_i\left| (\varphi _t^{*}J_M-J_{{\bf C}^n})(\frac \partial {\partial x^i})z^\alpha \right| _t\nonumber\\
&\leq&C_8\left| z\right| ^2(t+1)^{-\frac{2\varepsilon}{1+\varepsilon}}\ ,
\end{eqnarray}for some positive constant $C_8$ depending only on $\varepsilon,C_1,C_2$
and $n$.

For any fixed $t\geq 1,$ we define%
$$
r(t)=(\frac{\bar c}2)^{\frac 12}(t+1)^{\frac{\varepsilon}{2(1+\varepsilon)}} 
$$
and consider the $\bar \partial -$equation 
\begin{equation}
\label{3.25}\bar \partial \xi ^\alpha (z,t)=\bar \partial z^\alpha ,\qquad
z\in \hat B(0,r(t)),\quad \alpha =1,2,\cdots ,n\ . 
\end{equation}

By (\ref{3.24}), we know 
\begin{equation}
\label{3.26}\left| \bar \partial ^tz^\alpha \right| _t\leq C_8(\frac{\bar c}%
2)(t+1)^{-\frac \varepsilon{1+\varepsilon}}\ ,\qquad on\quad \hat B(0,r(t))\ ,\quad \alpha
=1,2,\cdots ,n\ . 
\end{equation}
When the positive constant $\bar c$ is chosen small enough, the estimates (%
\ref{3.11})-(\ref{3.16}) imply the metric $g_{ij}^{*}(\cdot ,t)$ can be
arbitrarily close to the standard K\"ahler metric on $B\left(0,\bar c(t+1)^{\frac
\varepsilon{1+\varepsilon}}\right).$ Thus the sphere $\partial \hat B(0,r(t))$ must be strictly convex
with respect to the metric $g_{ij}^{*}(\cdot ,t)$ (i.e. the second
fundamental form of $\partial \hat B(0,r(t))$ with respect to metric $%
g_{ij}^{*}(\cdot ,t)$ is strictly positive).

Then by using $L^2$ estimate theory for $\bar \partial -$operator ( see the
book of H\"ormander \cite{Ho} ), from (\ref{3.6}), (\ref{3.11})-(\ref{3.17}) and (%
\ref{3.26}) we know that (\ref{3.25}) has smooth solutions \{$\xi ^\alpha
(z,t)|\alpha =1,2,...n$\} such that for $\alpha =1,2,\cdots ,n,$%
\begin{equation}
\label{3.27}\left| \xi ^\alpha (z,t)\right| \leq \frac{C_9}{r(t)} 
\end{equation}
and 
\begin{equation}
\label{3.28}\left| \hat \nabla \xi ^\alpha (z,t)\right| \leq \frac{C_9}{%
r(t)^2} 
\end{equation}
on $\hat B(0,r(t)),$ where $C_9$ is a positive constant depending only on $%
\varepsilon,C_1,C_2$ and $n.$

We now define a new map \ $\Phi _t=(\Phi _t^1,\Phi _t^2,...\Phi _t^n):$ $%
\hat B(0,r(t))\rightarrow$ {\bf C}$^n$ by%
$$
\Phi _t^\alpha =z^\alpha -\xi ^\alpha (z,t)\ ,\qquad \alpha =1,2,\cdots ,n\
. 
$$

The equation (\ref{3.25}) says that $\Phi _t$ is holomorphic from $\hat
B(0,r(t))$\ equipped with the pull back complex structure $\varphi _t^{*}J_M$%
\ to {\bf C}$^n$ equipped with \ standard complex structure. When $t$ is large
enough , the estimate (\ref{3.27}) and (\ref{3.28}) imply $\Phi _t$ is a
diffeomorphism from $\hat B(0,r(t))$ to $\Phi _t(\hat B(0,r(t)))(\subset${\bf C}$^n)$ and 
\begin{equation}
\label{3.29}\hat B\left(0,\frac 12r(t)\right) \subset \Phi _t(\hat
B(0,r(t)))\subset \hat B(0,2r(t))\ . 
\end{equation}
Thus the map $\Phi _t\circ \varphi _t^{-1}$ is a holomorphic and injective
map from $B_t(x_0,r(t))(\subset M)$ to {\bf C}$^n$ (equipped with the standard
complex structure). It follows direcly from (\ref{3.29}) that the image of
the ball $B_t(x_0,r(t))$ of radius $r(t)$ with respect to the metric $%
g_{\alpha \bar \beta }(\cdot ,t)$ contains the Euclean ball $\hat B\left(0, 
\frac{r(t)}2\right) (\subset${\bf C}$^n).$

Denote%
$$
\Omega _t=\left( \Phi _t\circ \varphi _t^{-1}\right) ^{-1}\left( \hat
B\left(0,\frac{r(t)}2\right) \right) \ . 
$$
As noted before, the Ricci flow (\ref{2.1}) is shrinking, so by the virtue
of (\ref{3.29}) , for any $t$,%
$$
\Omega _t\supset B_t\left( x_0,\frac{r(t)}4\right) \supset B_0\left( x_0, 
\frac{r(t)}4\right) \ . 
$$

Hence there exists a sequence of $t_k\rightarrow +\infty ,$ as $k\rightarrow
\infty ,$ such that%
$$
M=\bigcup\limits_{k=1}^\infty \Omega _{t_k}\qquad and\qquad \Omega
_{t_1}\subset \ \ \Omega _{t_2}\subset \cdot \cdot \cdot \ . 
$$
Since each\ $\Omega _{t_k}(k=1,2,...)$ is biholomorphic to the unit ball of
{\bf C}$^n,$ $(\Omega_{t_k},\Omega_{t_l})$ is Runge pair for any $%
k,l. $ Then we can appeal to a theorem of Markoe \cite{Ma} (~see also Siu \cite{Si} )
to conclude that $M$ is a Stein manifold. This completes the proof of
Theorem1.1.

Note that $\Omega _{t_k}(k=1,2,...)$ are a\ sequence of exhaustion domains
of $M$ such that each of them is biholomorphic to the unit ball of {\bf C}$^n.$
Then by the Steiness of $M,$ one can repeat the gluing argument of Shi in
Section $9$ of \cite{Sh4} to construct a biholomorphic map from $M$ to a
pseudoconvex domain of {\bf C}$^n.$ In his paper \cite{Sh4}, Shi assumed the positivity
of Riemannian sectional curvature of $M$ to ensure the Steiness of the
manifold. Thus we get the proof for the part (2) of Theorem1.2 .

Therefore we have completed the proof of Theorem 1.1 and Theorem1.2.

\section*{4. A Gap Theorem}

\setcounter{section}{4} \setcounter{equation}{0}\ ~\quad In this section we
are interested in the question that how much the curvature could have near
the infinities for complete noncompact K\"ahler manifolds with nonnegative
bisectional curvature. The curvature behavior of a complete Riemannian
surface is quite arbitrary. For example, it is easy to construct complete
metrics on ${{\bf R}}^2$ from surface of revolution such that their curvatures
are zero outside some compact set, nonnegative everywhere, and positive
somewhere. It is surprising that the corresponding situation can not occur
for higher dimentions. The isometrically embedding theorem of Mok, Siu and
Yau stated before implies that there cannot too less positive bisectional
curvature near the infinity for a (complex) $n-$dimensional K\"ahler
manifold with $n\geq 2$. Recently, Ni \cite{N} extended the Mok-Siu-Yau's result
to some nonmaximal volume growth K\"ahler manifolds.
Here we give a further generalization as follows:\vskip  3mm

{\bf \underline{Theorem 4.1}} \ \ \ Suppose $M$ is a complete noncompact
K\"ahler manifold of complex dimension $n\geq 2$ with bounded and
nonnegative holomorphic bisectional curvature. Suppose there exists a
positive function $\varepsilon:{{\bf R}}\rightarrow {{\bf R}}$ with$\lim \limits_{r\rightarrow
+\infty }\varepsilon(r)=0$, such that for any $x_0,$%
$$
\frac 1{Vol(B(x_0,r))}\int\nolimits_{B(x_0,r)}R(x)dv\leq \frac{\varepsilon(r)}{r^2}\ . 
$$
Then $M^n$ is isometrically biholomorphic to a flat complete K\"ahler
manifold.\vskip  3mm

{\bf \underline{Proof.}} \ \ \ Suppose the metric in the theorem is $\tilde
g_{\alpha \bar \beta }(\cdot ).$ Consider the Ricci flow 
\begin{equation}
\label{4.1}\left\{ 
\begin{array}{ll}
\bigbreak \displaystyle\frac \partial {\partial t}g_{\alpha \overline{\beta }%
}(x,t)=-R_{\alpha \overline{\beta }}(x,t)\ , & x\in M\ ,\ t>0\ , \\ 
g_{\alpha \bar \beta }(x,0)=\tilde g_{\alpha \bar \beta }(x)\ , & x\in M\ . 
\end{array}
\right. 
\end{equation}
Let $F(x,t)$ be the nonpositive function defined in (\ref{2.2}). In Section\
2 we have known that $F(x,t)$ satisfies the following two differential
inequalities 
\begin{equation}
\label{4.2}e^{F(x,t)}\frac{\partial F(x,t)}{\partial t}\geq \Delta
_0F(x,t)-R(x,0)\ , 
\end{equation}
and 
\begin{equation}
\label{4.3}\Delta _0F(x,t)\leq R(x,0)\ ,\qquad x\in M\ ,\quad t\geq 0\ , 
\end{equation}
where $\triangle _0$ is the Laplace operator with respect to the initial
metric $\tilde g_{\alpha \bar \beta }$.

Exactly as in the proof of (\ref{2.18}), we can get from (\ref{4.2}) that
for any $x_0\in M$ and any $a>0,$%
\begin{equation}
\label{4.4}\frac 1{Vol(B_0(x_0,a))}\int\nolimits_{B_0(x_0,a)}F(x,t)dV_0\geq
- \frac{C_1t}{a^2}(1-F_{\min }(t))^2\ , 
\end{equation}
where $C_1$ is a positive constant depending only on $n$ and the function $%
\varepsilon(r).$

While by using (\ref{4.3}), a slight modification of the proof for (\ref
{2.15}) gives the estimate 
\begin{equation}
\label{4.5}F(x_0,t)\geq -\tilde \varepsilon(a)\log (2+a)+\frac{C_2}{Vol(B_0(x_0,a))}%
\int\nolimits_{B_0(x_0,a)}F(x,t)dV_0 
\end{equation}
for$\ $any $a>0,$ $x_0\in M$ and $t>0.$ Here $\tilde \varepsilon(r)$ is some
continuous positive function with $\lim \limits_{r\rightarrow +\infty
}\tilde \varepsilon(r)=0$ and $C_2$ is a positive constant depending only on $n$ and the
function $\varepsilon(r).$

Choose $x_0\in M$ such that $F(x_0,t)\leq \frac 12F_{\min }(t).$ We get from
(\ref{4.4}) and (\ref{4.5}) that for any $a>0,$%
\begin{equation}
\label{4.6}F_{\min }(t)\geq -2\tilde \varepsilon(a)\log (a+2)-\frac{2C_1C_2t}{a^2}%
(1-F_{\min }(t))^2\ . 
\end{equation}
By taking $a=(2+t)(1-F_{\min }(t)),$ we get 
\begin{equation}
\label{4.7}F_{\min }(t)+\log (1-F_{\min }(t))\geq -C_3\left[\tilde{\tilde\varepsilon}(t)\log (2+t)+1\right] 
\end{equation}
where $C_3$ is a positive constant depending only on $n$ and $\varepsilon(r),\ \tilde{\tilde\varepsilon}(t)$
is a positive function with $\lim \limits_{r\rightarrow +\infty}\tilde{\tilde\varepsilon}(r)=0$.

In particular, it follows from (\ref{4.7}) 
\begin{equation}
\label{4.8}\lim \limits_{t\rightarrow +\infty }\frac{F_{\min }(t)}{\log (2+t)%
}=0\ . 
\end{equation}

On the other hand , by the differential Harnack inequality of Cao \cite{Cao} as
before, we have%
\begin{eqnarray}
\label{4.9}
-F(x,2t)&=&\int_0^{2t}R(x,s)ds\nonumber\\
&\geq&tR(x,t)\int_t^{2t}\frac 1sds\nonumber\\
&=&(\log 2)\cdot tR(x,t)\ .
\end{eqnarray}The combination of (\ref{4.8}) and (\ref{4.9}) implies that
the solution $g_{\alpha \bar \beta }(\cdot ,t)$ of the Ricci flow (\ref{4.1}%
) exists for all time $t>0.$By using the differential Harnack inequality again, we get for $t>1,$%
\begin{eqnarray}
\label{4.10}
-F(x,t)&=&\int_0^{t}R(x,s)ds\nonumber\\
&\geq&\int_{\sqrt{t}}^tR(x,s)ds\nonumber\\
&\geq&\sqrt{t}R\left(x,\sqrt{t}\right)\int_{\sqrt{t}}^t\frac 1sds\nonumber\\
&=&\left(\frac 12\log t\right)\cdot\left(\sqrt{t}R\left(x,\sqrt{t}\right)\right)\ ,
\end{eqnarray}which together with (\ref{4.8}) implies 
\begin{equation}
\label{4.11}\lim \limits_{t\rightarrow \infty }\sqrt{t}R\left( \cdot ,\sqrt{t%
}\right) =0\ , 
\end{equation}

While $\sqrt{t}R\left( \cdot ,\sqrt{t}\right) $ is nondecreasing in time by
the differential Harnack inequality, therefore we conclude that%
$$
R(x,t)=0\qquad on\quad M\times [0,+\infty )\ . 
$$
This completes the proof of the theorem.\hfill
$Q. E. C. $

\section*{5. A Remark on Shi's Flat Metric}

\setcounter{section}{5} \setcounter{equation}{0}\ ~\quad In \cite{Sh3}, Shi
studied the following Ricci flow on K\"ahler manifolds, 
\begin{equation}
\label{5.1}\left\{ 
\begin{array}{ll}
\bigbreak \displaystyle\frac \partial {\partial t}g_{\alpha \overline{\beta }%
}(x,t)=-R_{\alpha \overline{\beta }}(x,t)\ , & x\in M\ ,\ t>0\ , \\ 
g_{\alpha \bar \beta }(x,0)=\tilde g_{\alpha \bar \beta }(x)\ , & x\in M\ . 
\end{array}
\right. 
\end{equation}
where $(M,\tilde g_{\alpha \bar \beta })$ is a complete noncompact K\"ahler
manifold of complex dimension $n$ with positive holomorphic bisectional
curvature.

He introduced the following Harnack quantity%
$$
Q(x,t)=\left\{ 1+g^{\alpha \bar \delta }(x,t)g^{\gamma \bar \beta
}(x,t)g^{\xi \bar \varsigma }(x,0)\tilde \nabla _\xi g_{\alpha \bar \beta
}\tilde\nabla _{\bar\zeta} g_{\gamma \bar \delta }(x,t)\right\} ^{\frac 12}\ , 
$$
where $\tilde \nabla $ is the covariant derivative with respect to the
initial metric $\tilde g_{\alpha \bar \beta }.$ By a direct computation, one
can get
\begin{equation}
\label{5.2}\frac{\partial Q}{\partial t}\leq C(n)\left| \nabla _\gamma
R_{\alpha \bar \beta }\right| _t\ ,
\end{equation}
where $\left| \cdot \right| _t$ is the norm with respect to the evolving
metric $g_{\alpha \bar \beta }(\cdot ,t).$ As seen in the previous sections,
under the assumptions that for a fixed $x_0\in M,$%
\begin{eqnarray}
(i) & Vol(B_0(x_0,r))\geq C_1r^{2n}\ ,& 0\leq r<+\infty\ ,\nonumber\\
(ii)&R(x,0)\leq\frac{C_2}{1+d(x_0,x)^2}\ ,
& x\in M\ ,\nonumber
\end{eqnarray}
the Ricci flow (\ref{5.1}) exists for all time $t\geq 0.$ Furthermore Shi \cite{Sh3}
proved the quantity $Q(x,t)$ is uniformly bounded on $M\times [0,+\infty ).$

Let $S(M)$ be the sphere bundle with respect to the initial metric $%
\tilde g_{\alpha \bar \beta }.$ Since the holomorphic bisectional curvature is
positive , the holonomy group is transitive on $S(M).$ For any two elements $%
(x,v),(y,w)\in S(M),$ suppose $\gamma :[0,1]\rightarrow M$ is a piecewise
smooth curve such that $\gamma (0)=x,\gamma (1)=y$ and $V(s)$ is a parallel
vector field along $\gamma $ such that $V(0)=v,V(1)=w.$ At the smooth point $%
\gamma (s)$ of $\gamma ,$ let%
$$
\theta (s)=\frac{\gamma ^{\prime }(s)}{\left| \gamma^{\prime }(s)\right| } 
$$
be the unit tangent vector of $\gamma $ at $\gamma(s).$

Since $Q$ is bounded , it follows from the definition of $Q$ that
\begin{eqnarray}
\left| \tilde \nabla _{\theta (s)}\log g_{V(s)\bar V(s)}(\gamma
(s),t)\right| _0^2&=&\frac 1{\ [g_{V(s)\bar V(s)}(\gamma (s),t)]^2}\left[
\tilde \nabla _{\theta (s)}g_{V(s)\bar V(s)}(\gamma (s),t)\right] ^2\nonumber\\
&\leq&const\ .\nonumber
\end{eqnarray}
Integrating from $\gamma (0)$ to $\gamma (1)$ along $\gamma $, we get
\begin{equation}
\label{5.3}\left| \log \frac{g_{w\bar w}(y,t)}{g_{v\bar v}(x,t)}\right| \leq
const\cdot \left| \gamma \right| \ ,
\end{equation}
where $\left| \gamma \right| $ is the length of $\gamma $ with respect to
the initial metric $\tilde g_{\alpha \bar \beta }.$

We now fix a point $(x_0,v_0)\in S(M)$ and define%
$$
U(t)=g_{v_0\bar v_0}(x_0,t)\ ,\qquad 0\leq t<+\infty  
$$
and
$$
\hat g_{\alpha \bar \beta }(x,t)=\frac 1{U(t)}\ \ g_{\alpha \bar \beta
}(x,t)\ ,\qquad on\quad M\times [0,\infty )\ . 
$$
Since Ricci flow (\ref{5.1}) is shrinking, we have $U(t)\leq 1$ for $t\in
[0,+\infty ).$ Thus the estimates (\ref{3.5}),(\ref{3.6}) on the curvature
of \ \ $\hat g_{\alpha \bar \beta }(x,t)$ and its derivatives still hold .
Combining with (\ref{5.3}) one can choose a sequence of times $%
t_k\rightarrow +\infty $ such that%
$$
\hat g_{\alpha \bar \beta }(\cdot ,t_k)\rightarrow\hat g_{\alpha \bar
\beta }(x,\infty )\ , 
$$
in the $C_{loc}^\infty $ topology of $M,$ where \ \ $\hat g_{\alpha \bar
\beta }(x,\infty )$ is a smooth flat metric on $M$ (by the decay estimate(%
\ref{3.5})). Clearly one cannot ensure the completeness for the limit metric 
$\hat g_{\alpha \bar \beta }(x,\infty ).$

From the above argument, one can see that once we have a bound on the
Harnack\ quantity $Q$, then we can construct a flat K\"ahler metric on the
manifold. \ However, in order to get a bound on $Q,$ the quadratic decay
assumption (ii) can be weakened to the following condition%
$$
(ii)^{\prime }\qquad R(x,0)\leq \frac{C_2}{1+d_0(x,x_0)^{\left( \frac
32+\varepsilon\right) }}\ ,\qquad x\in M\ , 
$$
where $\varepsilon$ is an arbitrarily small positive constant .

In fact by the estimate (\ref{3.5}) in Section3, we have%
\begin{equation}
\label{5.4}\left| \nabla _\gamma R_{\alpha \bar \beta }\right| _t\leq C_3(t+1)^{-\frac{%
3+6\varepsilon}{3+2\varepsilon}}\ ,\qquad on\quad M\times [0,+\infty )\ . 
\end{equation}
Since $Q(\cdot ,0)=0,$ we deduce from (\ref{5.2}) and (\ref{5.4}) that%
$$
Q(x,t)\leq \int_0^\infty C(n)\left| \nabla _\gamma R_{\alpha \bar \beta
}\right| _tdt\leq C_4\ , 
$$
where $C_4$ is some positive constant depending only on $\varepsilon,C_1,C_2$ and $n.$

\ Hence we obtain the following result.\vskip  3mm

{\bf \underline{Proposition 5.1}} \ \ \ Let $M$ be a complete noncompact
K\"ahler manifold of complex dimension $n$ with positive holomorphic
bisectional curvature. Suppose the maximal volume growth condition (i) \ and 
$\frac 32-$decay assumption (ii)$^{^{\prime }}$ hold for $M.$ Then with
respect to the complex structure we can construct a flat K\"ahler metric on $%
M$ . \ \

\end{document}